\documentclass{article}
  \newtheorem{proposition}{Proposition}[section]
  \newtheorem{theorem}[proposition]{Theorem}
  \newtheorem{lemma}[proposition]{Lemma}
  \newtheorem{corollary}[proposition]{Corollary}
  \newtheorem{definition}[proposition]{Definition}
  \newtheorem{example}[proposition]{Example}
  
  \newtheorem{remark}[proposition]{Remark}

  \newcommand{\mathbb}[1]{{\bf #1}}
\usepackage{eufrak}
  \begin{document}
  \title{Some applications of Andr\'e-Quillen homology to classes of arithmetic
rings}
  \author{Tiberiu Dumitrescu and Cristodor Ionescu\\[2mm]
  Facultatea de Matematica, Universitatea Bucuresti, \\ Str.
  Academiei 14, Bucharest, RO-010014, Romania\\
  email tiberiu@al.math.unibuc.ro\\[2mm]
  Institute of Mathematics \textit{Simion Stoilow} of the Romanian Academy,\\
  P.O. Box 1-764, Bucharest, RO 014700, Romania
  \\
  email Cristodor.Ionescu@imar.ro}
  \date{}
  \maketitle
  \begin{abstract}
We compute the first Andr\' e-Quillen homology modules for the simple
over-rings of integrally
closed domains and study an ideal theoretic condition arising from the
vanishing of $H_1$.
\end{abstract}
\vspace{5mm}

Andr\'e-Quillen (co)homology is known to be a powerful tool in
characterizing various
classes of rings or morphisms between noetherian rings.
Regular and complete intersection local rings, regular,
(formally) smooth or complete intersection morphisms can be
characterized with the help of this theory
(see Andr\'e (1974) and Brezuleanu et al. (1993)).
Classes of arithmetical integral domains, such as Pr\"ufer
domains, have also been characterized in this way (see Planas-Vilanova (1996)).

Let $A$ be an integrally closed
domain with quotient field $K$, $0\neq a,b\in A$ and $B=A[a/b]$.
In section $1$, we compute  $H_0(A,B,B)$ 
and $H_1(A,B,B),$ and we describe
$H_2(A,B,B)$ (Theorems \ref{2}  and \ref{15}). In particular, we show that
$H_1(A,B,B)=0$ if and only if $a^2A\cap b^2A=(aA\cap bA)^2$.

In section $2$, we investigate this condition in its own.
We say that $D$ is a  {\em $\star$-domain} if
$a^2D\cap b^2D=(aD\cap bD)^2$  for every $a,b\in D$.
The locally GCD domains are typical examples of $\star$-domains.
In Proposition \ref{998} we characterize the $\star$-pseudo-valuation domains.
In Corollary \ref{1113} we show that a  two-generated domain
(e.g. a quadratic extension of {\bf Z})  is a $\star$-domain if and only if it is Dedekind.
Finally, in Proposition \ref{72}, we prove that  the local class group of
a Krull  $\star$-domain has no element of order two.

Throughout this paper all rings are commutative.
For any undefined notation or terminology, the reader is refered
to Andr\'e (1974) and Gilmer (1972).

\section{Homological results}

In Theorem \ref{2} we compute the first two Andr\'e-Quillen
homology modules for a simple
over-ring of an integrally closed domain. By an over-ring of an integral domain $A$, we mean any intermediate ring between $A$ and its quotient field.  We need the following
well-known result, cf. Gilmer (1972, Corollary 34.9).
\begin{lemma}\label{1}
Let $A$ be an integrally closed domain with quotient field $K$
and let $0\neq f\in K[X]$. Then 
$$fK[X]\cap A[X]=fFA[X]$$
where $F=\{d\in K|\ df\in A[X]\}$.
\end{lemma}

\begin{theorem}\label{2}
Let $A$ be an integrally closed
domain with quotient field $K$, $0\neq a,b\in A$ and $B=A[a/b]$.
Let $I$ be the kernel of the
$A$-algebra morphism $\pi:R\rightarrow B$
sending $X$ to $a/b$, where $R=A[X]$.
We identify $B$ with $R/I$.
Then
  $$(i)\ \ \ \ \ \ \ \ \Omega_{B/A}\cong \frac{B}{(bA:_Aa)B} \cong
\frac{abA}{(aA+bA)(aA\cap bA)}
  \otimes_AR$$
  $$(ii)\ \ \ \ \ \ \ \ \ \ \ \ \ \ \ H_1(A,B,B)\cong \frac{a^2A\cap
b^2A}{(aA\cap bA)^2}\otimes_AR.$$
  \end{theorem}
  {\bf Proof.} $(i).$
  Set $c=a/b$.
  Since $A$ is integrally closed, Lemma \ref{1} shows that
  \begin{equation}\label{3}
  I=(X-c)K[X]\cap R=(X-c)(A:_A c)R=(X-c)(bA:_Aa)R.
  \end{equation}
  In particular, $I=(X-c)R\cap R.$
  Since $\Omega_{R/A}\otimes_R B\cong B$, the Jacobi-Zariski sequence
induced by
  $A\hookrightarrow R\stackrel{\pi}\rightarrow B$ is
  \begin{equation}\label{10}
  0\rightarrow H_1(A,B,B)\rightarrow
  I/I^2\stackrel{\delta}{\rightarrow}
  B\rightarrow \Omega_{B/A}\rightarrow 0.\end{equation}
  Let $q\in I$. Then $q=(X-c)r$ for some $r\in R$.
  Denoting the derivative of a polynomial $h\in R$ by $h'$,
  we get
  $$b\delta(q+I^2)=\delta(bq+I^2)=(bq)'(c)=br(c)$$ so
  \begin{equation}\label{8}
  \delta(q+I^2)=r(c) \mbox{ where } r=\frac{q}{X-c}.
  \end{equation}
  Hence
  $$\mbox{Im}(\delta)=(bA:_Aa)B.$$
  So
  $$\Omega_{B/A}\cong B/\mbox{Im}(\delta)\cong
  \frac{B}{(bA:_Aa)B}
  \cong
  \frac{R}{(X-c)(bA:_Aa)R+(bA:_Aa)R}\cong $$
  $$
  \cong \frac{R}{(R+Rc)(bA:_Aa)}\cong
  \frac{abR}{(aR+bR)(aR\cap bR)}
  \cong
  \frac{abA}{(aA+bA)(aA\cap bA)}\otimes _AR.$$
  $(ii).$
  Now let $f\in I$. By $(\ref{8})$, $\delta(f+I^2)=0$ if and only if
  $f\in (X-c)I$.
  So, by $(\ref{10})$,
  \begin{equation}\label{9}
  H_1(A,B,B)\cong\mbox{ker}(\delta)\cong \frac{(X-c)I\cap I}{I^2}.
  \end{equation}
  By $(\ref{3})$,
  we get
  \begin{equation}\label{11}
  (X-c)I\cap I=(X-c)^2K[X]\cap (X-c)R\cap R=(X-c)^2K[X]\cap R.
  \end{equation}
  Using Lemma \ref{1}, we get
  \begin{equation}\label{12}
  (X-c)^2K[X]\cap R=(X-c)^2(A[X]:_K (X-c)^2)R=(X-c)^2(b^2A:_Aa^2)R.
  \end{equation}
  Indeed, $A:_A c^2 \subseteq A:_A c$, because $A$ is integrally
closed. So
  $$A[X]:_K (X-c)^2=A:_K(1,2c,c^2)=A:_Ac^2=b^2A:_Aa^2.$$
  Combining
  $(\ref{9})$, $(\ref{11})$, $(\ref{12})$ and $(\ref{3})$,
  and taking account that $aA\cap bA=a(bA:_Aa)$, we get
  $$
  H_1(A,B,B)\cong \frac{(X-c)^2(b^2A:_Aa^2)R}{(X-c)^2(bA:_Aa)^2R}\cong
  \frac{a^2A\cap b^2A}{(aA\cap bA)^2}\otimes_AR.\ \bullet$$

  \begin{remark}{\em
  a) Note that $\Omega_{B/A}=B/(bA:_Aa)B$ is a flat $B$-module only when it
is zero, that is, $(bA:_Aa)B=B$.

   b) The fact that the $R$-modules
  $$\frac{abA}{(aA+bA)(aA\cap bA)} \otimes_AR\ \mbox{ and }\ 
  \frac{a^2A\cap b^2A}{(aA\cap bA)^2}\otimes_AR$$
  are $B$-modules, that is, they are annihilated by $I,$ can be proved
  directly as follows.
  Note that $I\subseteq (bA:_Aa)R +(aA:_Ab)R$.
  By symmetry, it suffices to see that both modules are annihilated by
$bA:_Aa$.
  For the first module this is clear because
  $(bA:_Aa)ab\subseteq (aA+bA)(aA\cap bA)$.
  As $A$ is integrally closed, $a^2A\cap b^2A\subseteq a^2A\cap abA$.
  Hence $(a^2A\cap b^2A)(bA:_Aa)\subseteq (a^2A\cap abA)(bA:_Aa)=(aA\cap
bA)^2,$ so the second module
  is annihilated by $I$.
  }\end{remark}

  \begin{corollary}\label{24}
  Let $A$ be an integrally closed
  domain, 
  $0\neq a,b\in A$ and $B=A[a/b]$.
  \par\noindent i) The following conditions are equivalent:
  \par a) $\Omega_{B/A}=0;$
  \par b) $aA+bA$ is an invertible ideal of $A;$
  \par c) $B=(bA:_Aa)B;$
  \par d) $H_n(A,B,E)=0$ for every  $B$-module $E$
and for every $n\geq 0$.
  \par\noindent ii) Consequently, the following conditions are equivalent:
  \par a) $A$ is a Pr\" ufer domain;
  \par b) $\Omega_{A[y]/A}=0$ for every $y\in K;$
  \par c) $H_n(A,C,E)=0$ for every over-ring $C$ of $A$, every $C$-module $E$
and any $n\geq 0$.
  \end{corollary}
  {\bf Proof:}
  $i).$
  The equivalence of $a)$, $b)$ and $c)$
  follows directly from Theorem \ref{2} (see also Smith (1979, Theorem 1)).
  The implication $d)\Rightarrow a)$ is obvious.
  To complete, we prove that $b)$ implies $d)$.
  Assume that $aA+bA$ is an invertible ideal of $A$ and let $B$ be an
over-ring of $A$,
  $E$ an $B$-module and $n\geq 0$. We claim that $H_n(A,B,E)=0$.
  Localizing (Andr\'e, 1974, Corollaries 4.59 and 5.27),
  we may assume that $A$ is a quasi-local domain.
  It follows that $aA+bA$ is a principal, so it is generated by $a$ or
$b$.
  Hence $B$ is a localization of $A$, so $H_n(A,B,E)=0$
  cf. Andr\'e (1974, Corollary 5.25).

   $ii).$
  A domain $A$ is a Pr\" ufer domain if and only if
  every nonzero two-generated ideal of $A$ is invertible
  (cf. Gilmer (1972, Theorem 22.1)), so
  $ii)$ follows from $i)$.
  $\bullet$\\[2mm]

  Somewhat in the same vein, note that
  in Planas-Vilanova (1996) it was shown that $D$ is a Pr\"ufer domain if and only if
  $H_2(D,D/I,D/I)=0$ for each (three-generated) ideal $I$ of $D$.

  \begin{corollary}\label{25} If $A$ is an integrally closed
  domain, 
  $0\neq a,b\in A$ and $B=A[a/b]$, then
  $H_1(A,B,B)=0$ if and only if
  $a^2A\cap b^2A=(aA\cap bA)^2$.
  In particular, if $aA\cap bA$ is a flat ideal, then $H_1(A,B,B)=0.$
  \end{corollary}
  {\bf Proof:}
  The first assertion follows directly from Theorem \ref{2}.
  To complete, assume that
  the ideal $J=aA\cap bA$ is $A$-flat. Then
  $$J^2=(aA\cap bA)J=aJ\cap bJ=a^2\cap abA\cap b^2A.$$
  Since $A$ is integrally closed,
  $$a^2A\cap abA\cap b^2A=a^2A\cap b^2A,$$ because from
  $x\in a^2\cap b^2A$ it follows that $(x/ab)^2\in A$.
  Thus $J^2=a^2\cap b^2A$, so the first part applies.
  $\bullet$\\[2mm]

  \begin{remark}\label{900}
  {\rm a) Let $D$ be an integrally closed domain such that
$$(aD\cap bD)(cD\cap dD)=acD\cap adD\cap bcD\cap bdD, \  \forall a,b,c,d\in D$$
(see Zafrullah (1987)). The preceding proof shows that $H_1(D,E,E)=0$ for each
simple over-ring $E$ of $D$.

  \par\noindent b) By Corollary \ref{25}, we see that
  if $a,b,c,d\in A$ are nonzero elements
  of the integrally closed
  domain $A$ and $A[a/b]=A[c/d]$,
  then $a^2A\cap b^2A=(aA\cap bA)^2$ if and only if
  $c^2A\cap d^2A=(cA\cap dA)^2$.
  }
  \end{remark}
  \begin{remark}{\em
  Theorem \ref{2} is no longer true if $A$ is not integrally closed.
  Indeed, if $A=\mathbb{Z}[i\sqrt{3}]$ and $B=A[(1+i\sqrt{3})/2]$,
  then it is easy to see that $\Omega_{B/A}=0$
  but the ideal $(1+i\sqrt{3}, 2)A$ is not
  invertible.
  }\end{remark}

  \begin{example}{\em
  Consider the integrally closed domains
  $A=\mathbb{Z}[X]$ and $B=\mathbb{Z}[X/2]=A[X/2].$
  By Theorem \ref{2} and Corollary \ref{25}, 
  $$\Omega_{B/A}\cong
B/(2A:_AX)B\cong B/2B$$
 and
  $H_1(A,B,B)=0$.
  }\end{example}

  \begin{example}{\em
  \label{14}
  Let $K$ be a field and $A=K+xL[x]$ where $L=K(y)$ and $x,y$ are
indeterminates
  over $K$. By Anderson et al. (1991, Theorem 2.7)
   $A$ is integrally closed.
  Consider the ring $B=A[yx/x]=K[y]+xL[x]$.
  An easy computation shows that 
  $$yxA\cap xA=x^2L[x]$$   and
  $$(yx)^2A\cap x^2A=x^3L[x].$$
  By Theorem \ref{2} we obtain
  $$H_1(A,B,B)\cong \frac{(yx)^2A\cap x^2A}{(yxA\cap
xA)^2}\otimes_AB\cong
  \frac{x^3L[x]}{x^4L[x]}\otimes_AB\cong L\otimes_A B\cong L\otimes_K
K[y]\cong L[z] $$
  where $z$ is an indeterminate over $L$ and the $B$-module structure
on $L$ is given by the ring morphism
  $B\rightarrow L$ sending $f(x,y)$ into $f(0,y)$.
  Hence $H_1(A,B,B)$ is not finitely generated as a $B$-module.
  Similarly, we get 
  $$\Omega_{B/A}\cong B/xL[x]\cong K[y].$$
  }\end{example}
  %
  %
  %
  %
  \begin{remark}{\em Let $A$ be a
  domain, 
  $a,b\in A$ nonzero elements and let $B=A[a/b]$. Then
  $H^0(A,B,B)=\mbox{Der}_A(B,B)=0$. Indeed, if $D\in \mbox{Der}_A(B,B)$,
  then
  $$0=D(a)=D(b(a/b))=bD(a/b),$$ so $D(a/b)=0$, whence $D=0$.
   }\end{remark}

  \begin{theorem}\label{30}
  Let $A$ be an integrally closed
  domain, 
  $0\neq a,b\in A$ and $B=A[a/b]$.
  If $a^2A\cap b^2A=(aA\cap bA)^2$, then
  $H^1(A,B,B)=((bA:_Aa)B)^{-1}/B$.
  \end{theorem}
  {\bf Proof.} As $a^2A\cap b^2A=(aA\cap bA)^2$, Theorem \ref{2} gives
that $H_1(A,B,B)=0$.
  Then $H^1(A,B,B)=\mbox{Ext}^1_B(\Omega_{B/A},B)$, cf. Andr\'e (1974, Lemma 3.19).
  By Theorem \ref{2},
  we have the exact sequence
  $$ 0\rightarrow (bA:_Aa)B \rightarrow B \rightarrow \Omega_{B/A}
  \rightarrow 0$$
  which gives the exact sequence
  $$0\rightarrow \mbox{Hom}_B(B,B)\rightarrow
\mbox{Hom}_B((bA:_Aa)B,B)\rightarrow
  \mbox{Ext}^1_B(\Omega_{B/A},B)\rightarrow 0.$$
  Since $\mbox{Hom}_B(B,B)=B$ and
$\mbox{Hom}_B((bA:_Aa)B,B)=((bA:_Aa)B)^{-1}$,
  we have
  $\mbox{Ext}^1_B(\Omega_{B/A},B)=((bA:_Aa)B)^{-1}/B$.
  $\bullet$\\[2mm]

  For the next result we use the following notation.
  Let $C\subseteq D$ be an extension of domains, and $J$ an ideal of $D$.
  Consider the canonical $C$-module epimorphism
  $\pi_J:S_2( _CJ)\rightarrow J^2$ where
  $S_2( _CJ)$ denotes the degree-two part of the symmetric algebra of
$J$.
  Let $\alpha$ be a nonzero element of $D$. It is not hard to see
  $\ker (\pi_J)$ is $C$-isomorphic to $\ker (\pi_{\alpha J})$.
  Call the ideal $J$ {\em $C$-syzygetic} (Planas-Vilanova, 1996), if $\pi_J$ is an
isomorphism.
  \begin{theorem}\label{15}
  Consider the setup in Theorem \ref{2}.
  Then $$H_2(A,B,B)\cong W\otimes_AR$$
  where $W$ is the kernel of the canonical
  map $S_2( _A aA\cap bA)\rightarrow (aA\cap bA)^2$.
  In particular, $H_2(A,B,B)=0$ if and only if $aA\cap bA$ is a
  syzygetic ideal of $A$.
  \end{theorem}
  {\bf Proof.}
  We use the notations of the proof of Theorem \ref{2}.
  In the Jacobi-Zariski sequence induced by
  $A\hookrightarrow R\stackrel{\pi}\rightarrow B$
  $$
  H_2(A,R,B) \rightarrow H_2(A,B,B)\rightarrow H_2(R,B,B) \rightarrow
H_1(A,R,B)
  $$
  the extreme terms are zero, so $H_2(A,B,B)\cong H_2(R,B,B)$.
  Since $B=R/I$, we have $H_2(R,B,B)\cong \mbox{ker}(\pi_I)$,
  where $\pi_I$ is the canonical
  map $S_2(I)\rightarrow I^2$.
  By the proof of Theorem \ref{2}, $I=(X-c)(bA:_Aa)R$.
  As noted in the paragraph preceding
  Theorem \ref{15}, $\pi_I$ is $R$-isomorphic to
  $W\otimes_AR$
  where $W$ is the kernel of the canonical
  map $S_2( _A aA\cap bA)\rightarrow (aA\cap bA)^2$. $\bullet$

  \begin{example}{\em
  Consider the setup of Example \ref{14}. We have $yxA\cap xA=x^2L[x]$.
  By Theorem \ref{15} and the paragraph preceding,
  it follows that $$H_2(A,B,B)\cong W\otimes_AR$$ where $W$ is the kernel of
the canonical
  map $S_2( _AL[x])\rightarrow L[x]$.
  }\end{example}
  \begin{corollary}\label{16}
  In the setup in Theorem \ref{2} assume that $aA\cap bA=abA$ and let E be any B-module.
  Then:
  \par\noindent a) $\Omega_{B/A}\cong B/bB$;
  \par\noindent b) $H_1(A,B,E)\cong 0:_Eb$;
  \par\noindent c) $H_2(A,B,E)=0$;
  \par\noindent d) $H^1(A,B,E)\cong E/bE$;
  \par\noindent e) $H^2(A,B,E)=0$.

  \end{corollary}
  {\bf Proof.}
  a) Follows at once from Theorem \ref{2}.
  \par\noindent b),c) From Theorems \ref{2} and \ref{15}, it follows that $H_1(A,B,B)=H_2(A,B,B)=0.$ 
  Now by Andr\'e (1974, Lemma 3.19) we have
   $$H_i(A,B,E)\cong {\rm Tor}_i^B(\Omega_{B/A},E)\cong{\rm Tor}^B_i(B/bB,E), i=1,2.$$
    From the exact sequence
  $$0\rightarrow bB\rightarrow B\rightarrow B/bB\rightarrow 0$$
   we obtain the exact sequence
   $$0\rightarrow {\rm Tor}_1^B(B/bB,E)\rightarrow bB\otimes E\rightarrow B\otimes E\rightarrow B/bB\otimes E\rightarrow 0$$
   and this gives b). As for c), it follows from the fact that ${\rm fd}_B(\Omega_{B/A})\leq 1.$
  \par\noindent d) Again by Andr\'e (1974, Lemma 3.19),  we have that $$H^i(A,B,E)\cong {\rm Ext}^i_B(\Omega_{B/A},E),\ i=1,2.$$ Now everything folows from the exact sequence
  $$0\rightarrow{\rm Hom}(B/bB,E)\rightarrow{\rm Hom}(B,E)\rightarrow{\rm Hom}(bB,E)\rightarrow{\rm Ext}^1_B(B/bB,E)\rightarrow 0.$$
  Note that for $E=A$ the assertion follows also from Theorem \ref{30}.
  \par\noindent e) We have ${\rm pd}(\Omega_{B/A})\leq 1.$
  $\bullet$\\[2mm]

  \begin{remark}\label{555}
   {\rm If A is a GCD domain and B is a simple over-ring of A, then  it follows that there exists some element $b \in B$ such that the same formulas as in the preceding corollary hold.}
   \end{remark}

\section{Ideal theoretic results}

In this section, we consider the condition in Corollary \ref{25} in its own.
We give the following ad-hoc definition.
\begin{definition}\label{400}
Let $D$ be a domain with quotient field $K$.
Call a domain $D$ a  {\em $\star$-domain} if
$a^2D\cap b^2D=(aD\cap bD)^2$ (equivalently, $(bD:_Da)^2=b^2D:_Da^2$) for every $a,b\in D$.
\end{definition}

\begin{remark}\label{401}
{\rm a) The condition of being a $\star$-domain is
clearly local. Hence, a locally GCD domain is a $\star$-domain.

b) A $\star$-domain is $2$-root closed, that is,  $D$ contains every element
$x\in K$ such that $x^2 \in D$.
%
Indeed, if $0\neq a,b\in D$ such that
$(a/b)^2\in D$, then $D=b^2D:_Da^2=(bD:_Da)^2$, so $a/b\in D$.
In particular,
 $A[X^2,X^3]$ is not a
$\star$-domain, for any domain $A$.

c) While it is easy to see that ${\bf Z}[i\sqrt{3}]$ is $2$-root closed,
${\bf Z}[i\sqrt{3}]$ is not a $\star$-domain.
Indeed, let $\mathfrak{m}$ be the maximal ideal $(2, 1+i\sqrt{3})$ and set $E={\bf Z}[i\sqrt{3}]$.
Then $\mathfrak{m}(1+i\sqrt{3})\subseteq (2)$, so $(2):_E(1+i\sqrt{3})=\mathfrak{m}$.
Similarly, we get 
$$(2)^2:_E(1+i\sqrt{3})^2=(2):_E(1-i\sqrt{3})=\mathfrak{m}.$$
Hence 
$$(2)^2:_E(1+i\sqrt{3})^2\neq ((2):_E(1+i\sqrt{3}))^2.$$
See also Corollary \ref{1113} for a more general assertion.

d) A $2$-root closed domain $D$ such that $aD\cap bD$ is a flat ideal
for each $a,b\in D$ is a $\star$-domain.
Indeed, let $0\neq a,b\in D$. Since $D$ is a $2$-root closed domain,
it follows that $abD\subseteq a^2\cap b^2D$.
Now we may repeat the argument given in the proof of Corollary  \ref{25}.
}\end{remark}

According to Anderson and Dobbs (1980),
a quasi-local domain $(D,\mathfrak{m})$ with quotient field $K$ is
called a {\em pseudo-valuation domain} (PVD), if
$x,y \in K$ and $xy \in \mathfrak{m}$ imply $x \in \mathfrak{m}$ or $y \in \mathfrak{m}.$
Clearly, a valuation domain is a PVD.
Next, we  characterize the $\star$ PVDs
(for other equivalent assertions see Zafrullah (1987, Theorem 4.5)).

\begin{proposition}\label{998} Let $(D,\mathfrak{m})$ be a  PVD which is
not a valuation domain.
Then $D$ is a $\star$-domain
if and only if  $D$ is $2$-root closed and $\mathfrak{m}=\mathfrak{m}^2.$
\end{proposition}
{\bf Proof:} By part $b)$ of Remark \ref{401}, we may suppose that $D$ is $2$-root closed.
Let $aD$ and $bD$ be two incomparable principal ideals of $D.$  Then $a^2D$ and $b^2D$ are also incomparable (by the 2-root closedness). By Anderson and Dobbs (1980, Prop. 1.4), we get $bD:_Da=b^2D:_Da^2=\mathfrak{m}.$ Hence $D$ is a $\star$-domain if and only if $(bD:_Da)^2=b^2D:_Da^2$ for any incomparable  ideals $aD$ and $bD,$ if and only if $\mathfrak{m}=\mathfrak{m}^2.\bullet$\\[2mm]



Now let $(V,\mathfrak{m})$ be a valuation domain with residue field $L, k$ a proper
subfield
of $L$ and $D$ the pre-image of $k$ in $V.$ By Hedstrom and Houston (1978, Proposition 2.6), $D$ is a
PVD which is not a valuation domain. Applying the preceding proposition,
we get that D is a $\star$-domain  if and only if  $\mathfrak{m}=\mathfrak{m}^2$ and
$k$ is $2$-root closed in $L$
(i.e. $x \in L$ and $x^2 \in k \Rightarrow x \in k).$

According to Sally and Vasconcelos (1974), an ideal $I$
of a quasi-local domain $(D,\mathfrak{m})$
is said to be {\em stable} if $I^2=aI$ for some $a \in I.$

\begin{remark}\label{411}
{\rm
Let $(D,\mathfrak{m})$ be a quasi-local domain and $I$ an ideal of $D$.
As shown in the proof of Sally and Vasconcelos (1974, Theorem 3.4),
if $I$ and $I^2$ are generated by two elements, then $I$ is stable.
For the convenience of the reader we repeat the proof here.

Let $a,b$ generate $I.$
As $I^2=(a^2,ab,b^2)$ is two-generated, one of these 3 generators is
supefluous. If $I^2=(a^2,ab),$
  then $I^2=aI$, while if $I^2=(b^2,ab),$ then $I^2=bI.$
Assume that $ab \in (a^2,b^2)$ and $I^2$ is not equal to
$aI$ or $bI.$ Then $ab=ra^2+sb^2$ with $r,s \in \mathfrak{m}.$
Changing the pair $(a,b)$ to $(a-b,b)$, we get $I^2=aI.$}
\end{remark}

\begin{lemma}\label{1111} Let $(D,\mathfrak{m})$ be a quasi-local $\star$-domain.
If $\mathfrak{m}$ is stable, then $\mathfrak{m}$ is principal. 
\end{lemma}
{\bf Proof:}
Assume that $\mathfrak{m}$ is not principal.
As $\mathfrak{m}$ is stable, $\mathfrak{m}^2=a\mathfrak{m}$ for some $a \in \mathfrak{m}.$ Take $b \in \mathfrak{m},b\notin aD.$ Then
             $$\mathfrak{m}b \subseteq  \mathfrak{m}^2=a\mathfrak{m} \subseteq aD.$$
  Hence $aD:_Db=\mathfrak{m},$ because $b\notin aD.$ From $\mathfrak{m}^2=a\mathfrak{m}$ we get $\mathfrak{m}^3=a^2\mathfrak{m}.$
Then
             $$\mathfrak{m}b^2 \subseteq  \mathfrak{m}^3=a^2\mathfrak{m} \subseteq a^2D.$$
  Hence $a^2D:_Db^2=\mathfrak{m},$ because $D$ being a $2$-root closed domain implies
  $b^2 \notin a^2D.$
Since D is a $\star$-domain, we get
             $$\mathfrak{m}=a^2D:_Db^2=(aD:_Db)^2=\mathfrak{m}^2\subseteq aD.$$
So  $\mathfrak{m}=aD$, a contradiction. $\bullet$

\begin{proposition}\label{1112} Let $D$ be a locally Noetherian domain such that for each
maximal ideal $\mathfrak{m}$,
$\mathfrak{m}D_\mathfrak{m}$ and $(\mathfrak{m}D_\mathfrak{m})^2$ are generated by two elements.
Then $D$ is a $\star$-domain if and only if $D$ is almost Dedekind.
\end{proposition}

{\bf Proof}: Apply the lemma and the paragraph preceding it, and use the
fact that a local domain with
  principal maximal ideal is a DVR. $\bullet$

\begin{corollary}\label{1113} Let $A$ be  domain whose ideals are two-generated
(e.g. a quadratic extension of {\bf Z}). Then $D$ is a $\star$-domain if and only if $D$ is Dedekind.
\end{corollary}

We close by giving two results concerning Krull  $\star$-domains.

\begin{proposition}\label{71} A Krull domain  $D$ is a $\star$-domain
if and only if the square of every divisorial ideal of $D$ is also divisorial.
\end{proposition}

{\bf Proof}:
Assume that $D$ is a Krull domain and
let $X^{1}(D)$ be the set of all height-one primes of $D$. Denote
the divisorial closure of a nonzero ideal $I$ by $I_v$.
Let $0\neq a,b\in D$.
By Fossum (1973, Proposition 5.9),
$$((aD\cap bD)^2)_v =
\bigcap_{\mathfrak{p}\in X^{1}(D)}(aD\cap bD)^2D_\mathfrak{p}=
\bigcap_{\mathfrak{p}\in X^{1}(D)}(a^2D_\mathfrak{p}\cap b^2D_\mathfrak{p})=
a^2D\cap b^2D$$
because each $D_\mathfrak{p}$ is a DVR.
So $D$ is a $\star$-domain if and only if
$(aD\cap bD)^2$ is a divisorial ideal for each
$0\neq a,b\in D$.

Now let $I$ be an arbitrary nonzero divisorial ideal.
By Fossum (1973, Proposition 5.11), $I^{-1}=(c,d)_v$ for some $c,d\in I^{-1}$.
So $$I=I_v=(I^{-1})^{-1}=(c,d)^{-1}=c^{-1}D\cap d^{-1}D=p^{-1}(aD\cap bD)$$ for some
$a,b,p\in D\setminus\{0\}$.
Hence $I^2$ is divisorial if and only if $(aD\cap bD)^2$ is  divisorial.
The assertion follows.
$\bullet$
\\[2mm]

Let $D$ be a Krull domain and $Div(D)$ its group of divisorial fractional
ideals under the $v$-multiplication (Fossum, 1973, Proposition 3.4).
The {\em local class group} $G(D)$ of $D$ is the factor  group
$Div(D)$ modulo the subgroup of invertible fractional ideals (Bouvier and Zafrullah, 1988).
By Fossum (1973, Corollary 18.15), a Krull domain has zero local class group if and only if
it is locally factorial (hence it is a $\star$-domain).

\begin{corollary}\label{72} If $D$ is a Krull  $\star$-domain,
then $G(D)$ has no element of order two.
\end{corollary}

{\bf Proof}:
Let $I$ be a divisorial ideal of $D$ such that $(I^2)_v$ is invertible.
By the preceding proposition $I^2=(I^2)_v$, so $I$ is invertible.
$\bullet$

  \vspace{0.4cm} 
  {\bf Acknowledgement.} Work supported
  by the CEEX program of the Romanian Ministry of Education and
  Research, contract Cex05-D11-11/2005.


\begin{center}
{\bf REFERENCES}
\end{center}

\begin{list}{
 Anderson, D.D.,}\item Anderson, D.F., Zafrullah, M. (1991). 
Rings between $D[X]$ and $K[X]$, \textit{Houston J. Math.}
17:109-128.
\end{list}

\begin{list}{ Anderson, D.F.,}\item Dobbs, D.E. (1980). Pairs of rings with the same spectrum, 
\textit{Canad J. Math.} 32:362-384.
\end{list}

\begin{list}{ Andr\'e, M.}\item (1974). Homologie des alg\`ebres commutatives,
\textit{Springer, Berlin.}
\end{list}

 \begin{list}{ Bouvier, A.,}\item  Zafrullah, M. (1988).
On some class group of an integral domain, \textit{Bull. Soc. Math. Greece}
29:45-59.
\end{list}

\begin{list}{ Brezuleanu, A.,}  \item Radu, N., Dumitrescu, T. (1993).
Local Algebras, \textit{Ed. Universit\u a\c tii Bucure\c sti, Bucure\c
sti.}
\end{list}

\begin{list}{ Fossum, R.}\item (1973).
The Ideal Class Group of a Krull domain,
\textit{Springer,  New York}.
\end{list}

\begin{list}{ Gilmer, R.}\item (1972).
Multiplicative Ideal Theory,
\textit{Marcel Dekker, New York}.
\end{list}

\begin{list}{ Hedstrom,}\item J.R., Houston, E.G. (1978). Pseudo-valuation Domains, \textit{Pacific
J. Math.} 75:137-147.
\end{list}

\begin{list}{ Planas-Vilanova, }\item F. (1996). Rings of weak dimension one
and syzygetic ideals, \textit{Proc. Amer. Math. Soc.} 124:3015-3017.
\end{list}

\begin{list}{  Smith, }\item W.W. (1979).
Invertible ideals and overrings,
\textit{Houston J. Math.} 5:141-153.
\end{list}

\begin{list}{ Sally,}\item J., Vasconcelos, W. (1974). Stable rings,
\textit{J. Pure Appl. Algebra} 4:319-336.
\end{list}

\begin{list}{ Zafrullah,}\item M. (1987). On a property of pre-Schreier domains,
\textit{Comm. Algebra}  15:1895-1920.
\end{list}


  \end{document}